\documentclass[12pt]{amsart}

\input{xy}
\xyoption{all}
\usepackage{graphicx}
\usepackage{color}
\usepackage{verbatim}

\newtheorem{theorem}{Theorem} [section]
\newtheorem{prop}[theorem]{Proposition}

\newtheorem{conjecture}[theorem]{Conjecture}
\newtheorem{question}[theorem]{Question}
\newtheorem*{theoremA}{Theorem~A}

\theoremstyle{definition}

\theoremstyle{remark}
\newtheorem{remark}[theorem]{Remark}
\newtheorem*{remarks}{Remarks}

\numberwithin{equation}{section}
\numberwithin{figure}{section}

\newcommand\C{{\mathbb C}}
\newcommand\Chat { {\hat{\C}} } 

\renewcommand\P{{\mathbb P}}

\newcommand\Z{{\mathbb Z}}
\newcommand\Q{{\mathbb Q}}

\newcommand\cM{\mathcal{M}}

\newcommand\cB{\mathcal{B}}
\newcommand\cS{\mathcal{S}}

\renewcommand\phi{\varphi}
 
\newcommand\Gal{\operatorname{Gal}}

\newcommand\PSL{\mathrm{PSL}}

\newcommand\Qbar{\overline{\mathbb{Q}}}
\newcommand\kbar{\bar{k}}
\newcommand\Kbar{\overline{K}}

\newcommand\del{\partial} 
 
\newcommand\iso{\simeq} 

\newcommand\supp{\operatorname{supp}}   

\renewcommand\Re {\operatorname{Re}}

\renewcommand\Im {\operatorname{Im}} 

\newcommand\ord{\operatorname{ord}}  

\newcommand\Rat  {\mathrm{Rat}} 
\newcommand\Poly {\mathrm{Poly}} 
\newcommand\M {\mathrm{M}}

\newcommand\MPoly {\mathrm{MPoly}}

\newcommand\Per {\mathrm{Per}}

\begin{document}

\title{KAWA 2015 -- Dynamical moduli spaces and elliptic curves}

\author{Laura De Marco}

\email{demarco@math.northwestern.edu}

\date{\today}

\begin{abstract}
In these notes, we present a connection between the complex dynamics of a family of rational functions $f_t: \P^1\to \P^1$, parameterized by $t$ in a Riemann surface $X$, and the arithmetic dynamics of $f_t$ on rational points $\P^1(k)$ where $k = \C(X)$ or $\Qbar(X)$.  An explicit relation between stability and canonical height is explained, with a proof that contains a piece of the Mordell-Weil theorem for elliptic curves over function fields.  Our main goal is to pose some questions and conjectures about these families, guided by the principle of ``unlikely intersections" from arithmetic geometry, as in \cite{Zannier:book}.  We also include a proof that the hyperbolic postcritically-finite maps are Zariski dense in the moduli space $\M_d$ of rational maps of any given degree $d>1$.  These notes are based on four lectures at KAWA 2015, in Pisa, Italy, designed for an audience specializing in complex analysis, expanding upon the main results of \cite{BD:polyPCF, D:stableheight, DWY:Lattes}.

\bigskip
\noindent
\textsc{R\'esum\'e.}  Dans ces notes, nous donnons un lien entre la dynamique complexe d'une famille de fractions rationnelles $f_t: \P^1\to \P^1$, param\'etr\'ee par une surface de Riemann $X$, et la dynamique arithm\'etique de $f_t$ sur les points rationnels de $\P^1(k)$, o\`u $k = \C(X)$.  Une relation explicite entre stabilit\'e et hauteur canonique est \'etablie, avec une preuve qui contient une partie du th\'eor\`eme de Mordell-Weil pour les courbes elliptiques sur un corps de fonctions.  Notre but principal est de poser quelques questions et conjectures, guid\'es par le principe des ``unlikely intersections" en g\'eom\'etrie arithm\'etique (cf.~\cite{Zannier:book}).  Nous incluons aussi une preuve du fait que les applications hyperboliques postcritiquement-finies sont Zariski denses dans l'espace des modules $\M_d$ des applications rationnelles de degr\'e donn\'e $d > 1$. Ces notes sont bas\'ees sur un cours de 4 s\'eances donn\'ees \`a KAWA 2015 \`a Pise, Italie, destin\'ees \`a une audience sp\'ecialis\'ee en analyse complexe, et d\'eveloppent les principaux r\'esultats de \cite{BD:polyPCF, D:stableheight, DWY:Lattes}.
\end{abstract}



\maketitle

\thispagestyle{empty}

These notes are based on a series of four lectures at KAWA 2015, the sixth annual school and workshop in complex analysis, this year held at the Centro di Ricerca Matematica Ennio de Giorgi in Pisa, Italy.  The goal of the lectures was to explain some of the background and context for my recent research, concentrating on the main results of \cite{BD:polyPCF, D:stableheight, DWY:Lattes}, developing connections between one-dimensional complex dynamics and the arithmetic of elliptic curves. The text below is essentially a transcription of the lectures, with each section corresponding to one lecture, with a few added details and references.  

Lecture 1 introduces a connection between complex dynamics on $\P^1$ (specifically, the study of preperiodic points) and the Mordell-Weil theorem for elliptic curves (specifically, the study of torsion points).  In Lecture 2, I review fundamental concepts from complex dynamics, some of which are used already in Lecture 1.  I also present a question about stability and bifurcations.  Height functions are introduced in Lecture 3, and I sketch the proof of the main result from \cite{D:stableheight}.  Lecture 4 is devoted to a conjecture about the moduli space $\M_d$ of rational maps $f: \P^1\to \P^1$ of degree $d$, in the spirit of ``unlikely intersections" (compare \cite{Zannier:book}).  Finally, I have written an Appendix containing a proof that the hyperbolic postcritically finite maps are Zariski dense in the moduli space $\M_d$.  This fact is mentioned in Lecture 4, and while similar statements have appeared in the literature, I decided it would be useful to see a proof based on the ideas presented here.

The theme of these lectures is closely related to that of a lecture series by Joseph Silverman in 2010, and it is inspired by some of the questions discussed there; his notes were published as \cite{Silverman:moduli}.  Good introductory references include Milnor's book on one-dimensional complex dynamics \cite{Milnor:dynamics} and Silverman's books on elliptic curves and arithmetic dynamics \cite{Silverman:elliptic, Silverman:advanced, Silverman:dynamics}.

I would like to thank the organizers of the workshop, Marco Abate, Jordi Marzo, Pascal Thomas, Ahmed Zeriahi, and especially Jasmin Raissy for her tireless efforts to keep all of us educated and entertained.  I thank Tan Lei for asking the question, one year earlier, that led to the writing of the Appendix.  I also thank Matt Baker, Dragos Ghioca, and Joe Silverman for introducing me to this line of research.  Finally, I thank the referees for their helpful suggestions and pointers to additional references.  

\bigskip
\section{Complex dynamics and the Mordell-Weil Theorem}

Our main goal is to study the dynamics of holomorphic maps 
	$$f: \P^1 \to \P^1,$$
where $\P^1 = \P^1_\C = \Chat = \C \cup \{\infty\}$ is the Riemann sphere.  We are interested in features that persist, or fail to persist, in one-parameter families
	$$f_t : \P^1 \to \P^1.$$
Here, we view $t$ as a complex parameter, lying in a Riemann surface $X$, and we should assume that the coefficients of $f_t$ are meromorphic functions of $t\in X$.  Key questions are related to stability and bifurcations, to be defined formally in Lecture \ref{bifurcation}.  Roughly, we aim to understand the effect of small (and large) perturbations in $t$ on the dynamical features of $f_t$.  

The natural equivalence relation on these maps is that of {\em conformal conjugacy}.  That is, $f\sim g$ if and only if there is a M\"obius transformation $A$ so that $f = AgA^{-1}$.  

There are two key examples to keep in mind throughout these lectures.  The first is
	$$f_t(z) = z^2 + t$$
with $t\in \C$, the well-studied family of quadratic polynomials.  (Any degree 2 polynomial is conformally conjugate to one of this form.  Any degree 2 rational function with a fixed critical point is also conjugate to one of this form.)  The second example that will play an important role is
\begin{equation} \label{Lattes family}
	f_t(z) = \frac{(z^2-t)^2}{4z(z-1)(z-t)},
\end{equation}
a family of degree-4 Latt\`es maps, parameterized by $t \in \C\setminus\{0,1\}$.  We will say more about this example in \S\ref{Lattes computation}.

Note that any such family $f_t$, where the coefficients are meromorphic functions of $t$ in a compact Riemann surface $X$, may also be viewed as a single rational function 
	$$f \in k(z)$$
where $k$ is the field $\C(X)$ of meromorphic functions on $X$.  It is convenient to go back and forth between complex-analytic language and algebraic language.  For example, we will frequently identify a rational point $P = [P_1: P_2] \in \P^1(k)$, where $P_1$ and $P_2$ are not both zero and lie in the field $k$, with the holomorphic map 
	$$P:  X\to \P^1(\C)$$
defined by $P(t) = [P_1(t):P_2(t)]$.

\subsection{Torsion points on elliptic curves}  Background on complex tori and elliptic curves can be found in \cite{Silverman:elliptic} and \cite[Chapter 7]{Ahlfors}.  Recall that a number field is a finite extension of the rationals $\Q$.  An elliptic curve $E$, defined over a number field $k$, may be presented as the set of complex solutions $(x,y)$ to an equation of the form $\{y^2 = x^3 + A x + B\}$ with $A,B \in k$ satisfying $4A^3 + 27B^2 \not=0$.  (In fact, we should view $E$ as a subset of $\P^2$.)  By definition, the {\em rational points} $E(k)$ are the pairs $(x,y)\in k^2$ satisfying the given equation, together with the point $[0:1:0]\in\P^2$ representing the origin of the group.

In the 1920s, Mordell and Weil proved the following result:

\begin{theorem}  [Mordell, Weil, 1920s]  \label{MW}
If $E$ is an elliptic curve defined over a number field $k$, then the group of rational points $E(k)$ is finitely generated.
\end{theorem}

\noindent Mordell proved the theorem for $k=\Q$ and Weil generalized the result to number fields.

In particular, the theorem implies that {\em the set of torsion points in $E(k)$ is finite}.  Recall that a point $P\in E$ is torsion if there exists an $n$ so that $n\cdot P = 0$ in the additive group law on $E$.  In its usual complex presentation, we can write 
	$$E \iso \C/\Lambda$$
for a lattice $\Lambda$ in the complex plane.  The additive group law in $\C$ descends to the group law on $E$, and the torsion points of $\C/\Lambda$ are the rational combinations of the generators of the lattice.  We see immediately that there are infinitely many (in fact a dense set of) torsion points in $E$, when working over the complex numbers.

When $k$ is a function field, Theorem \ref{MW} as stated is false.  As an example, take $k = \C(t)$, the field of rational functions in $t$, and choose an $E$ such as $\{y^2 = x(x-1)(x-2)\}$, which has no dependence on $t$.  Over $\C$, the elliptic curve $E$ has infinitely many torsion points.  But all points constant in $t$ are rational in $k$ (that is, $E(\C) \subset E(k)$) so $E(k)$ contains infinitely many torsion points, and the group $E(k)$ is infinitely generated.  It turns out this kind of counterexample is the only way to create problems:

\begin{theorem} [Lang-N\'eron, 1959, Tate, N\'eron 1960s] \label{MW function field}
Let $E$ be an elliptic curve defined over a function field $k = \C(X)$ for a compact Riemann surface $X$, and assume that $E$ is not isotrivial.  Then $E(k)$ is finitely generated.  
\end{theorem}

To explain isotriviality, it helps to view an elliptic curve $E$ over $k = \C(X)$ as a family of complex curves $E_t$, $t\in X$, or indeed as a complex surface ${\bf E}$ equipped with a projection ${\bf E} \to X$, where the general fiber is a smooth elliptic curve.  The elements of $E(k)$ are simply the holomorphic sections of this projection.  The elliptic curve $E$ is {\em isotrivial} if the (smooth) fibers $E_t$ are isomorphic as complex tori.  In algebraic language, it will mean that after a base change, which corresponds to passing to a branched cover $Y\to X$ and pulling our fibers back to define a new surface ${\bf F} \to Y$, the new surface is birational to a product $Y\times E_0$.  In fact, Theorem \ref{MW function field} holds under the weaker hypothesis that $E$ is not {\em isotrivial over $k$}, meaning that the birational isomorphism to the product is over $k$, or with base $X$.  (We will revisit these notions of isotriviality in \S\ref{isotriviality}.)

For the proof of Theorem \ref{MW function field}, much of the strategy for the proof of Theorem \ref{MW} goes through; see \cite{Silverman:elliptic, Silverman:advanced} or \cite[Appendix B]{Baker:functionfields}.  A new ingredient is required for bounding the torsion part.  

\subsection{Torsion points are preperiodic points}  \label{torsion}
Take any elliptic curve $E$ over any field $k$ (well, let's work in characteristic 0 for simplicity).  Then the identification of a point $P$ with its additive inverse $-P$ defines a degree 2 projection 
\begin{equation} \label{projection}
	\pi: E(\kbar) \to \P^1(\kbar).
\end{equation}
This $\pi$ is given by the Weierstrass $\wp$-function, when $k=\C$ and $E = \C/\Lambda$.  Now let $\phi$ be an endomorphism of $E$.  
For example, let's take
	$$\phi(P) = P+P = 2P.$$
Note that a point $P$ is torsion on $E$ {\em if and only if} it has finite orbit under iteration of $\phi$.  That is, the sequence of points 
	$$P, \, 2P, \, 4P, \, 8P, \ldots$$
must be finite.  Descending to $\P^1$, the endomorphism $\phi$ induces a rational function $f_\phi$ so that the diagram 
$$\xymatrix{ E \ar[d]_\pi \ar[r]^\phi & E\ar[d]^\pi \\   \P^1   \ar[r]^{f_\phi} & \P^1  }$$	
commutes.  The degree of $f_\phi$ coincides with that of $\phi$, which is 4 in this example.  And the projection of a point $P$ is preperiodic for $f_\phi$ {\em if and only if} $P$ is torsion on $E$.  Recall that a point $x$ is preperiodic for $f$ if its forward orbit $\{f^n(x): n\geq 0\}$ is finite.

Rational functions $f: \P^1\to \P^1$ that are quotients $f = f_\phi$ of an endomorphism $\phi: E\to E$ are called {\em Latt\`es maps}.  A classification and summary of their dynamical features is given in \cite{Milnor:Lattes}.  

\subsection{Exploiting the dynamical viewpoint: an example}  \label{Lattes computation}
Consider the Legendre family of elliptic curves,
	$$E_t = \{y^2 = x(x-1)(x-t)\}$$
for $t\in \C\setminus\{0,1\}$, defining an elliptic curve $E$ over $k = \C(t)$.  It also defines an elliptic surface ${\bf E} \to X$ with $X =\P^1(\C)$; the fibers are smooth for all $t \in \C\setminus\{0,1\}$.  The projection $\pi : E(k) \to\P^1(k)$ of \eqref{projection} is given by $(x,y)\mapsto x$.  Take endomorphism $\phi(P) = 2P$ on $E$.  The action of $\phi$ on the $x$-coordinate induces a rational function $f \in k(x)$; namely, 
	$$f_t(x) = \frac{(x^2-t)^2}{4x(x-1)(x-t)},$$
the example from (\ref{Lattes family}).  See, e.g., \cite[Chapter III]{Silverman:elliptic}.

Take a point $P\in E(k)$, so $P$ projects via (\ref{projection}) to a point $x_P\in\P^1(k)$, and consequently it defines a holomorphic map 
	$$x_P: X \to \P^1,$$
where we recall that $X = \P^1$, so $x_P$ is a rational function in $t$.  For any $x\in\P^1(k)$, we may define 
	$$x_n(t) = f_t^n(x(t)).$$
We can compute explicitly from the formulas that if $x(t) = 0$, $1$, $t$, or  $\infty$, then $x_1(t) \equiv \infty$, so that $x_n(t) \equiv \infty$ for all $n\geq 1$.  These four points are the $x$-coordinates of the 2-torsion points of the elliptic curve $E$ defined over $k$.  On the other hand, for every other $x \in \P^1(k)$, it turns out that the topological degree of $x_1: X\to \P^1$ is at least 2, and then 
	$$\deg x_n = 4^{n-1} \deg x_1$$
for all $n\geq 1$ \cite[Proposition 3.1]{DWY:Lattes}.  In particular, $\deg x_n \to \infty$ with $n$, so $x_n$ cannot be preperiodic.  We conclude that
$$ |\{ P \in E(k): P \mbox{ is torsion}\}| = |\{x \in \P^1(k): x \mbox{ is preperiodic for } f\}| = 4.$$

\begin{remarks}
In general, a rational point of $E$ will descend to a rational point on $\P^1$ (though rational points on $\P^1$ typically do not lift to rational points on $E$), in a two-to-one fashion except over the four branch points of the projection $E\to \P^1$.  The original computation that $E(k)_{\mathrm{tor}}$ consists of exactly the 2-torsion subgroup in the Legendre family was probably done around 1900; the same result under the weaker hypothesis that only the $x$-coordinates of points in $E$ are rational would have required additional machinery, though may still have been known in the early 20th century \cite{Silverman:pc}.  
\end{remarks}

\subsection{A proof of finiteness in general}  \label{finiteness}
Here is a complex-dynamic proof of the finiteness of the rational torsion points, for an non-isotrivial elliptic curve defined over a function field, a key piece of Theorem \ref{MW function field}.  The main idea is to exploit the correspondence between torsion and preperiodic, as explained in \S\ref{torsion}.

Let $X$ be a compact Riemann surface.  Let $E$ be an elliptic curve defined over $k =\C(X)$, and assume that $E$ is not isotrivial.  As a family, $E_t$ is a smooth complex torus for $t$ in a finitely-punctured Riemann surface $V \subset X$ (where the discriminant is non-zero and the $j$-invariant is finite); non-isotriviality means that the $E_t$ are not all conformally isomorphic.  

The endomorphism $\phi(P) = 2P$ on $E$ induces a rational function $f\in k(z)$ of degree 4, via the projection (\ref{projection}).  Equivalently, we obtain a family $f_t$ of well-defined rational functions of degree 4 for all $t$ in the punctured Riemann surface $V$.  Non-isotriviality of $E$ guarantees that not all $f_t$ are conformally conjugate.  

Every rational point $P\in E(k)$ defines a holomorphic section $P: X\to {\bf E}$ of the complex surface, and then, via the quotient (\ref{projection}) from each smooth fiber $E_t \to \P^1$, $P$ defines a holomorphic map 
	$$x_P: X \to \P^1(\C).$$
From the discussion of \S\ref{torsion}, the point $P$ is torsion in $E$ if and only if the sequence of functions
	$$\{t \mapsto f_t^n(x_P(t))\}_{n \geq 1}$$
is finite.  We will study the (topological) degrees of the functions 
	$$x_n(t) = f_t^n(x(t))$$
from $X$ to $\P^1$ for every $x \in \P^1(k)$.

\begin{prop} \label{finite degree}
Fix $k = \C(X)$.  Let $f \in k(z)$ be any rational function of degree $\geq 2$.  There exists a constant $D = D(f)$ so that every preperiodic point $x \in \P^1(k)$ has degree $\leq D$ as a map $X \to \P^1$.  
\end{prop}

Proposition \ref{finite degree} follows from properties of the Weil height on $\P^1$ for function fields; see \S\ref{canonical} below.  It was first observed in \cite{Call:Silverman}.  Without appealing to height theory, this proposition can also be proved with an intersection-theory computation in the complex surface $X\times\P^1$, but we will skip the proof. 

By passing to a finite branched cover $Y\to X$, we may choose coordinates on $\P^1$ so that 
	$$f_t(\{0,1,\infty\}) \subset \{0,1,\infty\}$$
for all $t$ in a finitely-punctured $W\subset Y$.  For example, if $f_t$ has three distinct fixed points for all $t\in V$, we might place them at $0,1,\infty$ by conjugating $f$ with a M\"obius transformation $M \in k(z)$.  Unfortunately, this cannot be done in general, and we will need to pass to a branched cover of $V$ and add further punctures, in order to holomorphically label the fixed points and keep them distinct.  Moreover, if $f_t$ has $< 3$ fixed points for all $t\in V$, we choose at least one fixed point to follow holomorphically and two of its (iterated) preimages to obtain a forward-invariant set.  Setting $\ell = \C(Y)$, we may now view our $f$ as an element of $\ell(z)$, or equivalently, as an algebraic family of degree 4 rational functions $f_t$ with $t\in W$, where $W$ is the complement of finitely many points in $Y$.  

Take any point $x = x_0 \in \P^1(\ell)$ which is preperiodic for $f$, and look at the sequence $x_n(t) = f_t^n(x(t))$ of maps $Y \to \P^1$, for $n\geq 1$.  From Proposition \ref{finite degree}, there is a positive integer $D = D(f)$ so that $\deg x_n\leq D$ for all $n$.  Set
	$$S_n = x_n^{-1}\{0,1,\infty\} \cap W \subset Y.$$
Because of the normalization for $f$, with the set $\{0,1,\infty\}$ forward invariant, we see that $S_n \subset S_{n+1}$ for all $n$.  If $x_n$ is constant equal to 0, 1, or $\infty$ for some $n$, then $S_n = W$.  The degree bound of $D$ implies that the set
	$$S_x = \bigcup_{n\leq N_x} S_n \subset W$$ 
with 
	$$N_x = \max\{n: S_n \not= W\}$$
has cardinality $\leq 3 D$.  Now, for each $n\geq 0$, $x_n$ is a meromorphic function on $Y$ of degree $\leq D$; if nonconstant, it takes all zeroes, poles, and 1s in the finite set $S_x\cup (Y\setminus W)$; therefore, the number of distinct {\em nonconstant} functions in the sequence $\{x_n\}$ is bounded in terms of $D$ and the number of punctures in $W$.  Moreover, there cannot be more than 9 consecutive and distinct {\em constant} iterates of $x_0$; otherwise by interpolation (as in \cite[Lemma 2.5]{D:stableheight}), the degree 4 rational functions $f_t$ would be independent of $t$.  But independence of $t$ means that $f$ is isotrivial, which we have assumed is not the case.  

We conclude that the orbit length of $x_0$ is bounded by a number depending on $f$, but not depending on the point $x_0$ itself.  In other words, there is a constant $N = N(f)$ so that for every preperiodic $x \in \P^1(\ell)$, there exists $m = m(x) < N$ so that
	$$f_t^N(x(t)) = f_t^m(x(t))$$
for all $t\in W$.  Therefore, every preperiodic point $x \in \P^1(\ell)$ is a solution to one of $N$ algebraic equations over $\ell$ of degree $4^N$; in particular, the set of all preperiodic points in $\P^1(k)$ is finite.

\bigskip
\section{Stability and bifurcations}
\label{bifurcation}

Now we take a step back and introduce some important concepts from complex dynamics.  A basic reference is \cite{Milnor:dynamics}.  References for stability and the bifurcation current are given below.  

\subsection{Basic definitions}   \label{definitions}
Let $\P^1 = \Chat$ denote the Riemann sphere.  
Let $f: \P^1 \to\P^1$ be holomorphic.  Then $f$ is a rational function; we may write it as 
	$$f(z) = \frac{P(z)}{Q(z)}$$
for polynomials $P, Q\in\C[z]$ with no common factor.  Then the topological degree of $f$ is given by 
	$$\deg f = \max\{\deg P, \deg Q\}.$$
A {\em holomorphic family of rational maps} is a holomorphic map 
	$$f: X\times \P^1 \to \P^1$$
where $X$ is any complex manifold.  We write $f_t$ for the map $f(t, \cdot)$ for each $t\in X$.  Since $f$ is holomorphic, it follows that $\deg f_t$ is constant in $t$.  

\begin{remark}
This definition of holomorphic family not the same as the notion of family I introduced in Lecture 1.  There, I required $X$ to be compact of dimension 1 and the coefficients of $f_t$ to be meromorphic in $t$.  Such a family defines a holomorphic $f: V\times\P^1\to \P^1$ with $V$ the complement of finitely many points in $X$. 
\end{remark}

As an example, we could take $f: \C\times \P^1\to \P^1$ defined by $f_t(z) = z^2 + t$ with $t\in \C$.  A non-example would be $f_t(z) = tz^2 + 1$ with $t\in \C$, even though the coefficients are holomorphic.  This $f$ is not holomorphic at $(t,z) = (0, \infty)$; note that the degree drops when $t=0$.

An important example is the family of {\em all} rational maps of a given degree $d\geq 2$. Namely, 
	$$\Rat_d = \{f: \P^1\to \P^1 \mbox{ of degree } d\}$$
which may be seen as an open subset of $\P^{2d+1}(\C)$, identifying $f$ with the list of coefficients of its numerator and denominator.  In fact, $\Rat_d$ is the complement in $\P^{2d+1}(\C)$ of the {\em resultant hypersurface} consisting of all pairs of polynomials with a common root; thus $\Rat_d$ is a smooth complex affine algebraic variety.  Similarly, we define 
	$$\Poly_d = \{f: \C\to \C \mbox{ of degree } d\} \iso \C^*\times\C^d.$$

For any rational function $f$, its Julia set is defined by 
	$$J(f) = \Chat \setminus \Omega(f),$$
where $\Omega(f)$ is the largest open set on which the iterates $\{f^n\}$ form a {\em normal family}.  Remember what this means:  given any sequence of iterates, there is a subsequence that converges uniformly on compact subsets of $\Omega$.  Also recall Montel's Theorem:  if a family of meromorphic functions defined on an open set $U\subset \Chat$ takes values in a triply-punctured sphere, then it must be a normal family.   By Montel's Theorem, we can see that $J(f)$ may also be expressed as the smallest closed set such that $f^{-1}(J) \subset J$ and $|J|>2$.  Another useful characterization of $J(f)$ is as
	$$J(f)  = \overline{\{\mbox{repelling periodic points}\}}$$
where a point $x \in \P^1$ is periodic if there exists an $n\geq 1$ so that $f^n(x) = x$ and it is repelling if, in addition, $|(f^n)'(x)| > 1$.

\subsection{Stability}  \label{stable}
Suppose $f: X\times\P^1\to \P^1$ is a holomorphic family of rational maps.  Suppose that a holomorphic map
	$$c: X\to \P^1$$
parametrizes a critical point of $f$.  That is $f_t'(c(t)) = 0$ for all $t$.  The pair $(f,c)$ is {\em stable on $X$} if 
	$$\{t \mapsto f_t^n(c(t))\}_{n\geq 0}$$
forms a normal family on $X$.  An important fact about stability (see, e.g., \cite[Chapter 4]{McMullen:CDR}):  if $X$ is connected, then the pair $(f, c)$ is stable on $X$ if and only if either 
\begin{enumerate}
\item	$c(t)$ is disjoint from $J(f_t)$ for all $t\in X$; or 
\item $c(t) \in J(f_t)$ for all $t \in X$.
\end{enumerate}
In fact, it was proved by Lyubich and Ma\~n\'e-Sad-Sullivan \cite{Lyubich:stability, Mane:Sad:Sullivan} that $(f, c)$ is stable on $X$ for {\em all} critical points $c$ if and only if $f_t|J(f_t)$ is structurally stable on $X$.  This means that, on each connected component $U$ of $X$, the restrictions $f_t|J(f_t)$ are topologically conjugate for all $t\in U$.  In particular, the Julia sets are homeomorphic, and the dynamics are ``the same".  

For any holomorphic family $f: X\times\P^1\to \P^1$, the {\em stable locus} $\cS(f,c)\subset X$ is defined to be the largest open set on which the sequence of functions $\{t \mapsto f_t^n(c(t))\}_{n\geq 0}$ is normal.  A priori, this set could be empty, but it follows from the characterizations of stability that $\mathcal{S}(f,c)$ is always open and dense in $X$ \cite{Mane:Sad:Sullivan}.  We define the {\em bifurcation locus of the pair} $(f,c)$ by 
	$$\mathcal{B}(f,c) = X \setminus \cS(f,c),$$ 
and the {\em bifurcation locus of} $f$ as 
	$$\cB(f) = \bigcup_{c: f'(c)=0} \cB(f,c).$$

\begin{remark}  
Given a holomorphic family $f: X\times\P^1\to \P^1$, it may require passing to a branched cover of $X$ to holomorphically follow each critical point and define $\cB(f,c)$, but the union $\cB(f)$ is well defined on the original parameter space $X$, by the symmetry in its definition.  It coincides with the set of parameters where we fail to have structural stability on the Julia sets.
\end{remark}

\subsection{Example}
\begin{figure} [h]
\includegraphics[width=2.844in]{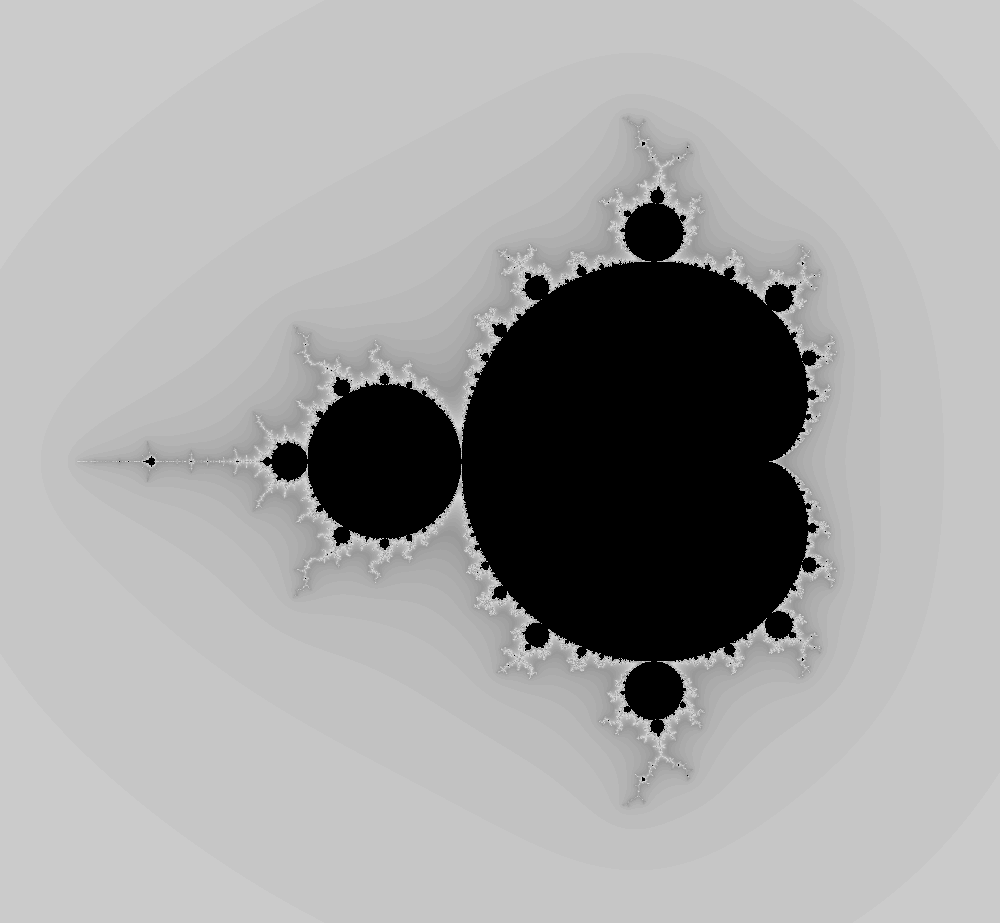}
\caption{ \small The Mandelbrot set.}
\label{Mandelbrot}
\end{figure}

Let $f_t(z) = z^2+t$ with $t\in \C$.  The bifurcation locus $\mathcal{B}(f)$ is equal to the boundary of the Mandelbrot set.  See Figure \ref{Mandelbrot}.  Indeed, the critical points of $f_t$ are at 0 and $\infty$ for all $t\in \C$.  Since $\infty$ is fixed for all $t$, it is clear that the family $\{t\mapsto f_t^n(\infty)\}_n$ is normal.   The Mandelbrot set is defined by 
	$$\cM = \{t \in \C:  \sup_n|f^n_t(0)| < \infty\}.$$
By Montel's Theorem, the pair $(f, 0)$ is stable on the interior of the Mandelbrot set, where the iterates are uniformly bounded.  It will also be stable on the exterior of the Mandelbrot set, because the iterates converge locally uniformly to infinity.  So any open set intersecting the boundary $\del \cM$ is precisely where the family $\{t\mapsto f_t^n(0)\}_n$ fails to be normal.  

\subsection{Bifurcation current}  \label{bifurcation current}
The bifurcation current was first introduced in \cite{D:current}; it quantifies the instability of a holomorphic family $f_t$.  Excellent surveys on the current and its properties are \cite{Dujardin:bifurcation} and \cite{Berteloot:lectures}.  

\begin{theorem}  \label{current theorem}
Let $f: X\times\P^1\to \P^1$ be a holomorphic family of rational maps of degree $d\geq 2$, and let $c: X\to \P^1$ be a critical point.  There exists a natural closed positive $(1,1)$-current $T_c$ on $X$, with continuous potentials, so that 
	$$\supp T_c = \cB(f,c).$$
\end{theorem}

The idea of its construction is as follows.  For each $t\in X$, there is a a unique probability measure $\mu_t$ on $\P^1_\C$ of maximal entropy for $f_t$, with support equal to $J(f_t)$.  Roughly speaking, the current $T_c$ is a pullback of the family of measures $\{\mu_t\}$ by the critical point $c$; as such, it charges parameters $t$ where the critical point $c(t)$ is ``passing through" the Julia set $J(f_t)$.  

For a family of polynomials $f_t$ of degree $d\geq 2$, we make this definition of $T_c$ precise by working with the escape-rate function
	$$G_t(z) = \lim_{n\to\infty} \frac{1}{d^n} \log^+ |f_t^n(z)|.$$
For each fixed $t$, the function $G_t$ is harmonic away from the Julia set, and it coincides with the Green's function on the unbounded component of $\C\setminus J(f_t)$ with logarithmic pole at $\infty$.  As such, its Laplacian (in the sense of distributions) is equal to harmonic measure on $J(f_t)$.  This is the measure $\mu_t$ of maximal entropy.  It turns out that $G_t(z)$ is continuous and plurisubharmonic as a function of $(t,z)$ on $X\times\C$ \cite{Branner:Hubbard:1}. The bifurcation current is defined on $X$ as
	$$T_c = (dd^c)_t \, G_t(c(t)).$$
When $X$ has complex dimension 1, the operator $dd^c$ is the Laplacian, and $T_c$ is a positive measure on $X$.  

For rational functions $f_t$, it is convenient to work in homogeneous coordinates on $\P^1$; then one can define an escape-rate function (locally in $t$) on $\C^2$ by 
\begin{equation} \label{escape rate}
	G_{F_t}(z,w) = \lim_{n\to\infty} \frac{1}{d^n} \log \| F_t^n(z,w) \|,
\end{equation}
where $F_t$ is a homogeneous presentation of $f_t$ and $\|\cdot \|$ is any choice of norm on $\C^2$.  Take any holomorphic (local) lift $\tilde{c}$ of the critical point into $\C^2\setminus\{(0,0)\}$, and set 
	$$T_c = (dd^c)_t \, G_{F_t}(\tilde{c}(t)).$$
The current is independent of the choices of $F$ and $\tilde{c}$, as long as they are holomorphic in $t$.

\subsection{Examples}
For the family $f_t(z) = z^2+t$ with $t\in \C$, the bifurcation current for $c(t)=0$ is (proportional to) the harmonic measure on the boundary of the Mandelbrot set.  For the family of Latt\`es maps, 
	$$f_t(z) = \frac{(z^2-t)^2}{4z(z-1)(z-t)},$$
with $t \in \C\setminus\{0,1\}$, the bifurcation measures are all equal to 0, because the critical points have finite orbit for all $t$; this family is stable on $X = \C\setminus\{0,1\}$.

\subsection{Questions} \label{bifurcation questions}
It has been an open and important problem, since these topics were first investigated, to understand and classify the stable components within the space $\Rat_d$ of all rational maps of degree $d$, or in $\Poly_d$, the space of all polynomials.  We still do not have a complete classification of stable components for the family $f_t(z) = z^2+t$, though conjecturally, all stable components will consist of {\em hyperbolic maps} (the quadratic polynomials for which there exists an attracting periodic cycle), and the hyperbolic components {\em have} been classified.  

In another direction, related to the topic of Lecture 4, we would like to know the answer to the following question:

\begin{question}  \label{critical question}
Suppose $X$ is a connected complex manifold and $f: X\times\P^1\to \P^1$ is a holomorphic family of rational maps, with marked critical points $c_1, c_2: X\to \P^1$.  Suppose the bifurcation loci $\cB(f,c_1)$ and $\cB(f,c_2)$ are nonempty, and assume that the bifurcation currents satisfy
	$$T_{c_1} = C\, T_{c_2}$$ 
as currents on $X$, for some constant $C>0$.  What can we conclude about the triple $(f, c_1, c_2)$?  
\end{question}

\noindent
As an example, if $c_1$ and $c_2$ share a grand orbit for all $t$, so that there exist integers $n, m \geq 0$ such that 
	$$f_t^n(c_1(t)) = f_t^m(c_2(t))$$
for all $t\in X$, then we may conclude that 
	$$T_{c_1} = d^{m-n} \, T_{c_2}$$
In general, we might expect that equality of the currents as in Question \ref{critical question} implies the grand orbits of $c_1$ and $c_2$ coincide, allowing for possible symmetries of $f$; see the discussion in \S\ref{unlikely} below and part (3) of Conjecture \ref{UI}.

Weakening the hypothesis of Question \ref{critical question}, we might ask about the bifurcation loci as sets:

\begin{question}  \label{bif question}
Suppose $X$ is a connected complex manifold and $f: X\times\P^1\to \P^1$ is a holomorphic family of rational maps, with marked critical points $c_1, c_2: X\to \P^1$.  Suppose the bifurcation loci $\cB(f,c_1)$ and $\cB(f,c_2)$ are nonempty, and assume that
	$$\cB(f,c_1) = \cB(f,c_2)$$ 
as subsets of $X$.  What can we conclude about the triple $(f, c_1, c_2)$?  
\end{question}

\noindent
As far as I am aware, there are no known examples where the bifurcation loci will coincide without orbit relations (in their general, symmetrized form) between $c_1$ and $c_2$.  In particular, equality of the bifurcation loci might imply equality of the currents $T_{c_1} = C \, T_{c_1}$ on $X$ for some constant $C>0$.  

\subsection{Bifurcation measures for arbitrary points} \label{point bifurcations}
Suppose $f: X\times\P^1\to \P^1$ is a holomorphic family of rational maps, and assume for simplicity that $\dim X = 1$.  In this case, the bifurcation currents are positive measures on $X$.    For arithmetic and algebraic applications, it is useful to have a notion of bifurcation currents for arbitrary points, not only the critical points.  The associated ``bifurcation locus" does not carry the same dynamical significance, in terms of topological conjugacy of nearby maps, but it does reflect the behavior of the point itself under iteration. 

\begin{figure} [h]
\includegraphics[width=3.5in]{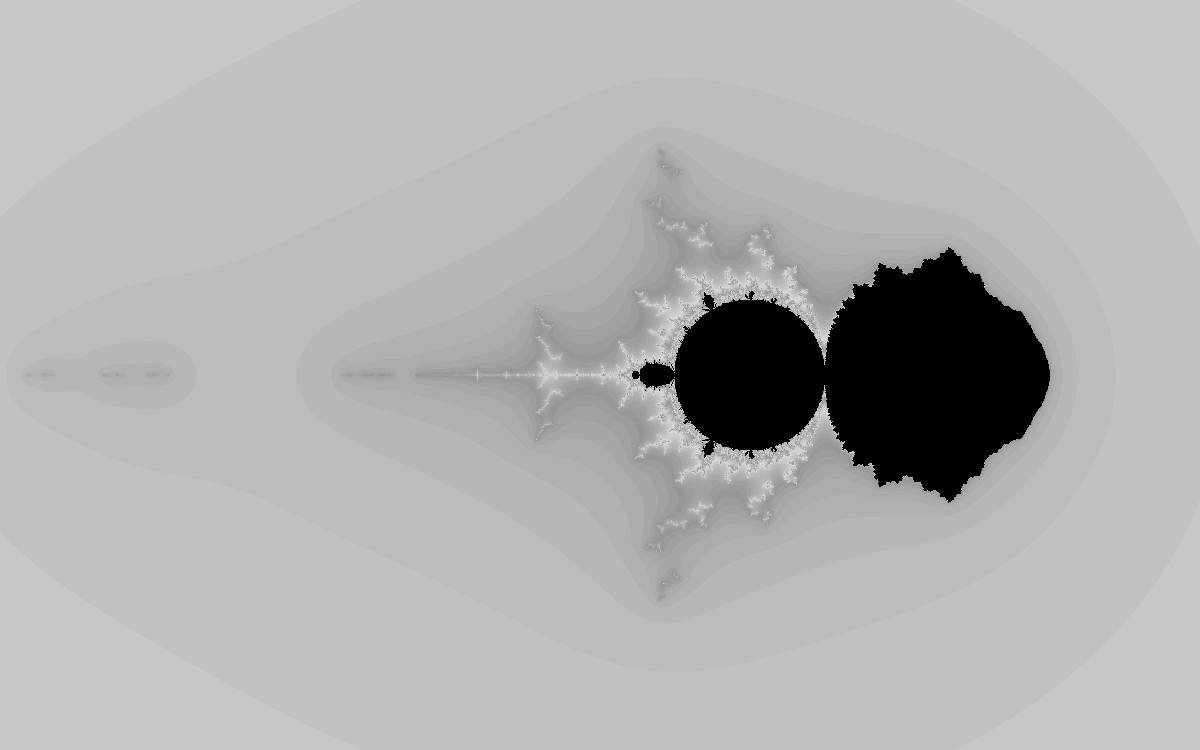}
\includegraphics[width=2.625in]{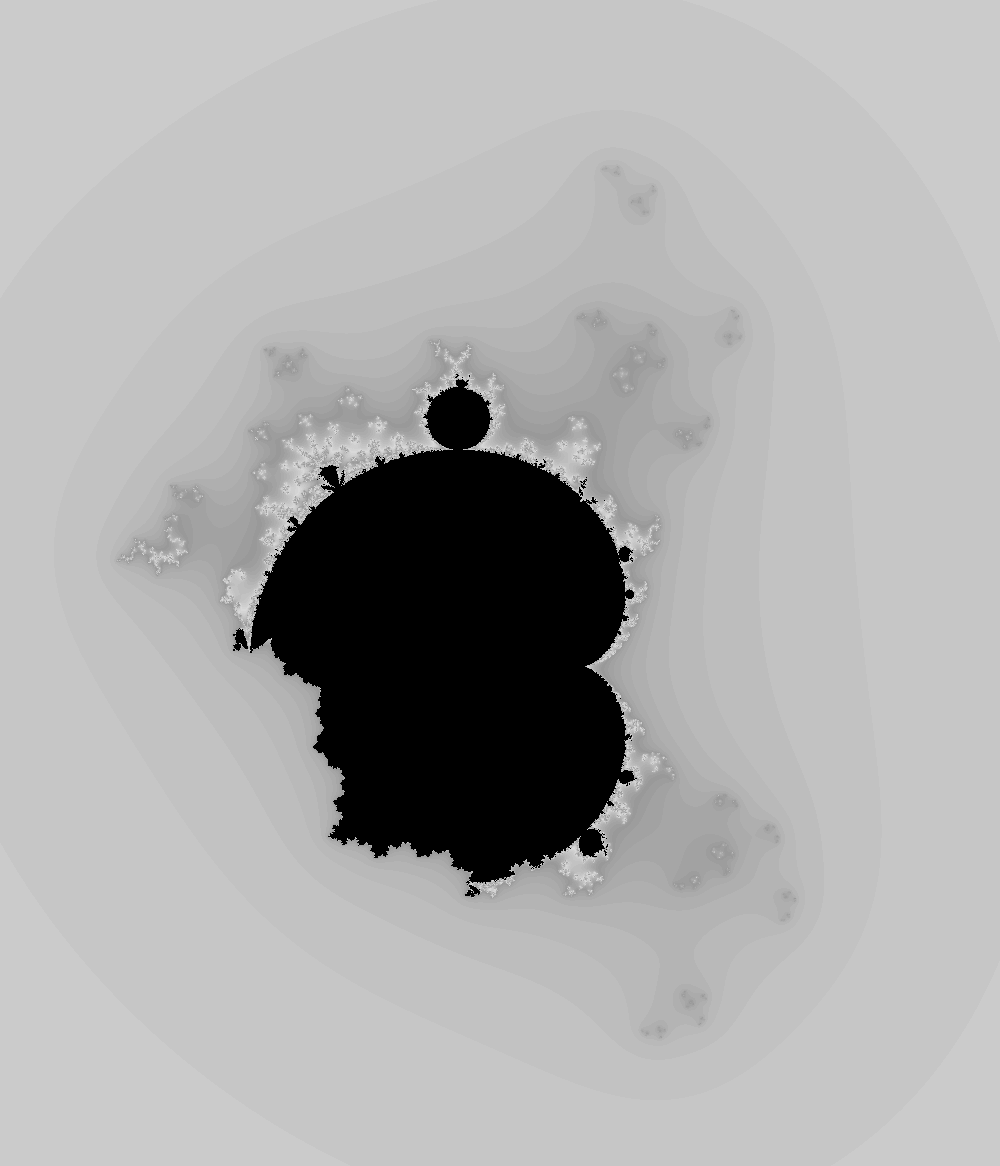}
\caption{ \small Let $f_t(z) = z^2 + t$ with $t\in \C$.  At left, with $P(t) = 1$, the bifurcation locus of $(f,P)$ is the boundary of the (disconnected) compact set shown here in the region $\{-3.5 \leq \Re t \leq 0.5, \; -1.25 \leq \Im t \leq 1.25\}$.  At right, with $P(t) = -0.122561 + 0.744862 i$ (which is an approximation of the center of the ``rabbit" hyperbolic component, a root of the equation $t^3 + 2t^2 + t + 1 = 0$ and a $t$ value where 0 has period 3), in the region $\{-1.5 \leq \Re t \leq 1.5, \; -1.5 \leq \Im t \leq 2.0\}$.  The two images here are shown at the same scale as Figure \ref{Mandelbrot}.  
}
\label{point 1}
\end{figure}

Let $P: X\to \P^1$ be any holomorphic map.  As in the case of critical points, we say that {\em the pair $(f, P)$ is stable on $X$} if the sequence 
	$$\{t \mapsto f_t^n(P(t))\}_{n\geq 1}$$
forms a normal family on $X$.  If it is not stable, then its stable locus $\cS(f,P)$ is the largest open set on which the sequence of iterates is normal; the bifurcation locus $\cB(f,P)$ is the complement, $X\setminus \cS(f,P)$.  The current $T_P$, which is a measure since $X$ has dimension 1, is locally expressed as $dd^c$ of the subharmonic function $G_{F_t}(\tilde{P}(t))$, exactly as described in \S\ref{bifurcation current} for critical points.  As before, we will have $\supp T_P  = \cB(f,P)$.  (A proof was given in \cite[Theorem 9.1]{D:lyap}.)  Two examples are shown in Figure \ref{point 1}.

\bigskip
\section{Canonical height}

Throughout this lecture, we let $X$ be a compact Riemann surface, and let $k = \C(X)$ be the field of meromorphic functions on $X$.  Suppose $f \in k(z)$ is a rational function of degree $d \geq 2$ with coefficients in $k$.  Away from finitely many points of $X$, say on $V = X \setminus\{x_1, \ldots, x_n\}$, $f$ determines a holomorphic family of rational maps $f: V\times \P^1\to\P^1$ of degree $d$.   In this setting, we will say that $f$ defines an {\em algebraic family} on $V$.  

\subsection{Isotriviality} \label{isotriviality}
We say $f \in k(z)$ is {\em isotrivial} if all $f_t$ for $t\in V$ are conformally conjugate.  Equivalently, there exists a M\"obius transformation $M\in \kbar(z)$ so that $M f M^{-1}$ is an element of $\C(z)$.  In other words, there is a branched cover $p: W \to V$ and an algebraic family of M\"obius transformations $M: W \times \P^1 \to \P^1$ so that 
	$$M_s f_{p(s)} M_s^{-1}: \P^1 \to \P^1$$
is independent of the parameter $s \in W$.    

We say that $f\in k(z)$ is {\em isotrivial over $k$} if there exists a degree 1 function $M \in k(z)$ so that $M f M^{-1} \in \C(z)$.  In other words, there is no need to pass to a branched cover $W \to V$ to define the family of M\"obius transformations $M: V\times \P^1 \to \P^1$.  

A function $f$ can be isotrivial but fail to be isotrivial over $k$ only if it has non-trivial automorphisms.  For example, consider the cubic polynomial $P(z) = z^3 - 3 z$.  This polynomial $P$ commutes with $z\mapsto -z$.  Conjugating by $M_s(z) = sz$ so that $M_s^{-1}PM_s(z) = s^2 z^3 - 3 z$, and setting $t = s^2$, we see that the family 
	$$f_t(z) = t z^3 - 3 z$$
with $t \in \C^*$ is isotrivial over a degree-2 extension of $k = \C(t)$.  But it is not isotrivial over $k$, because a conjugacy between $f_1$ and $f_t$ with $|t|=1$ induces a nontrivial monodromy as $t$ moves around the unit circle \cite[Example 2.2]{D:stableheight}.

\subsection{Stability and finite orbits}
In his thesis, McMullen showed the following; it was a key ingredient in his study of root-finding algorithms (such as Newton's method).  Recall the definitions of stability and bifurcations from the previous lecture.

\begin{theorem} \label{McMullen} \cite{McMullen:families}
Suppose $f\in k(z)$ is a rational function of degree $>1$, and assume that $f$ is not isotrivial.  If the bifurcation locus $\cB(f)\subset V$ is empty, then $f$ is a flexible family of Latt\`es maps; that is, there exists a non-isotrivial elliptic curve $E/\kbar$ and an endomorphism $\phi: E\to E$, so that the diagram 
$$\xymatrix{ E \ar[d]_\pi \ar[r]^\phi & E\ar[d]^\pi \\   \P^1   \ar[r]^{f} & \P^1  }$$	
commutes, where $\pi$ is the projection (\ref{projection}).
\end{theorem}

Theorem \ref{McMullen} is actually a statement about the critical orbits of $f$.  First, we pass to a branched cover $Y\to X$ if necessary so that the critical points $c_1, \ldots, c_{2d-2}$ can be labelled holomorphically along a punctured Riemann surface $W\subset Y$ where $f_t$ is well defined.  Then, we recall from Lecture 2 that $\cB(f) = \emptyset$ if and only if the sequence of maps 
	$$\{t\mapsto f_t^n(c_i(t))\}_{n\geq 1}$$
forms a normal family on $W$ for each $i$.  McMullen observed that this, in turn, is equivalent to the condition that each critical point must be persistently preperiodic along $V$.  That is, for each $i$, there exist positive integers $n_i>m_i$ so that 
	$$f_t^{n_i}(c_i(t)) \equiv f_t^{m_i}(c_i(t))$$
for all $t\in W$.  Finally, he combined this with Thurston's work on postcritically finite maps, who showed that the only non-isotrivial families of postcritically finite maps are Latt\`es families \cite{Douady:Hubbard:Thurston}.  

More recently, Dujardin and Favre extended McMullen's theorem to treat the critical points independently.  For their study of the distribution of postcritically finite maps in $\Rat_d$ and $\Poly_d$, they proved:

\begin{theorem} \label{DF critical}  \cite{Dujardin:Favre:critical}
Let $f: V\times\P^1 \to \P^1$ be any non-isotrivial algebraic family of rational maps with critical point $c: V\to \P^1$.  The pair $(f,c)$ is stable on $V$ if and only if $(f,c)$ is preperiodic.  
\end{theorem}

The proofs of Theorems \ref{McMullen} and \ref{DF critical} use crucially the fact that the point is critical.  But, it is not a necessary assumption for the conclusion of the theorem: 

\begin{theorem}  \cite{D:stableheight}  \label{mine}
Suppose $f\in k(z)$ has degree $d\geq 2$ and is not isotrivial over $k$.  Let $V = X\setminus\{x_1, \ldots, x_n\}$ be the punctured Riemann surface so that $f: V\times\P^1\to\P^1$ defines an algebraic family.  Fix any $P \in\P^1(k)$.  The following are equivalent.
\begin{enumerate}
\item	The pair $(f, P)$ is stable on $V\subset X$; i.e., the sequence $\{t\mapsto f_t^n(P(t))\}_{n\geq 1}$ forms a normal family on $V$.
\item The pair $(f, P)$ is preperiodic; i.e., there exist positive integers $n>m$ so that $f_t^n(P(t)) = f_t^m(P(t))$ for all $t$ in $V$.
\item The canonical height $\hat{h}_f(P)$ is equal to 0.
\end{enumerate}
Moreover, the set 
	$$\{P\in\P^1(k): (f,P) \mbox{ is stable}\}$$
is finite.  
\end{theorem}

The definition of the canonical height $\hat{h}_f$ will be given in \S\ref{canonical}.  First some remarks.

The equivalence (2) $\iff$ (3) and the finiteness of preperiodic rational points were first proved by Matt Baker \cite{Baker:functionfields}, under the slightly stronger hypothesis of non-isotriviality (versus the non-isotriviality over $k$ appearing here).   His proof uses the action of $f$ on the Berkovich analytification of $\P^1$, working over an algebraically closed and complete extension of the field $k$, with respect to each place of $k$; his methods extend to any product-formula field $k$, not only the geometric case considered here.

When $f$ is a Latt\`es family such as (\ref{Lattes family}), Theorem \ref{mine} contains a piece of the Mordell-Weil theorem, Theorem \ref{MW function field}, since the preperiodic points for $f$ coincide with the torsion points on the elliptic curve, as discussed in \S\ref{torsion}.

\subsection{The logarithmic Weil height}  
To introduce heights, it is instructive to begin with the logarithmic Weil height $h$ on $\P^1(\Qbar)$.  Suppose first that the field $k$ is $\Q$, the field of rational numbers.  We set
	$$h(x/y) = \log \max\{|x|, |y|\}$$
where $x/y$ is written in reduced form with $x,y\in \Z$.  Recall that a number field is a finite extension of $\Q$; we extend $h$ to number fields as follows:  for non-zero element $\alpha$ of a number field $k$, we consider the minimal polynomial $P_\alpha(x) = a_0 x^n + \cdots + a_n$ in $\Z[x]$ and we set 
	$$h(\alpha) = \frac{1}{n}\; \left( \log|a_0|  +  \sum_{\{\alpha': P_\alpha(\alpha')=0\}} \log\max\{|\alpha'|,1\} \right).$$
An alternative expression for $h$ is given in terms of the {\em set of places} $\cM_k$ of the field $k$: the places $v$ of a number field are an infinite collection of inequivalent absolute values on $k$ and come equipped with multiplicities $n_v$ so that they satisfy the {\em product formula}
	$$\prod_{v\in \cM_k} |\alpha|_v^{n_v} = 1$$
for all $\alpha\in k^*$.  For example, with $k=\Q$, there is a $p$-adic absolute value for every prime $p$ and the standard Euclidean absolute value, often denoted $|\cdot |_\infty$ to connote the ``infinite" prime; we have $n_v = 1$ for all places $v$ of $\Q$.  

The logarithmic Weil height on $\Qbar$ is equal to
	$$h(\alpha) = \frac{1}{[k: \Q]} \sum_{v\in \cM_k} n_v \log\max\{|\alpha|_v,1\}$$
for any number field $k$ that contains $\alpha \in \Qbar$.  This definition carries over to points $(x:y)\in \P^1(\Qbar)$ as 
	$$h(\alpha:\beta) = \frac{1}{[k: \Q]} \sum_{v\in \cM_k} n_v \log\max\{|\alpha|_v, |\beta|_v\}.$$
The sum is independent of the choice of representative $(x,y) \in k^2$ because of the product formula.

These expressions for the height $h$ make sense for other fields equipped with a product formula, such as the function field $k = \C(X)$ with $X$ a compact Riemann surface.  The {\em places} of $k = \C(X)$ now correspond to the complex points of the Riemann surface $X$, where 
	$$|\phi|_x = e^{-\ord_x\phi},$$
with $\ord_x$ computing the order of vanishing (or of a pole, when negative) at the point $x$.  The product formula 
	$$\prod_{x \in X} |\phi|_x = 1$$
holds for all non-zero meromorphic functions $\phi$, because the total number of zeroes must equal the total number of poles in $X$, counted with multiplicity.  In this case, we may take all multiplicities $n_x$ to be equal to 1.  We set	
	$$h(\alpha:\beta) = \sum_{x \in X} \log\max\{|\alpha|_x, |\beta|_x\}$$
for all $(\alpha:\beta) \in \P^1(k)$.  In fact, from the definition of $|\cdot|_x$ for each $x\in X$, the Weil height is nothing more than  
\begin{equation} \label{degree}	
	h(\phi) = \deg (\phi: X \to \P^1)
\end{equation}
for every $\phi \in \P^1(k)$, because the topological degree of $\phi$ is the total number of poles of $\phi$, counted with multiplicity.  

\subsection{The canonical height}  \label{canonical}
In \cite{Call:Silverman}, Call and Silverman introduced a canonical height function for polarized dynamical systems on projective varieties, defined over number fields or function fields, based on Tate's construction of the canonical height for elliptic curves.  The height of a point represents its arithmetic complexity; a canonical height should reflect the growth of arithmetic complexity along an orbit.  We present the definition only for a rational map $f: \P^1\to \P^1$ of degree $d\geq 2$.  

If $f: \P^1\to \P^1$ is a rational map of degree $d\geq 2$, defined over a field $k$ equipped with the product formula, the {\em canonical height} of $f$ is defined as 
	$$\hat{h}_f(P) = \lim_{n\to\infty} \frac{1}{d^n} h(f^n(P))$$
where $h$ is the logarithmic Weil height and $P \in \P^1(k)$.  The limit exists and the function $\hat{h}_f$ is uniquely determined by the following two properties \cite{Call:Silverman}:
\begin{enumerate}
\item[(H1)] $ \hat{h}_f(f(P)) = d \, \hat{h}_f(P)$ for all $P \in \P^1(k)$, and
\item[(H2)] There is a constant $C = C(f)$ so that $|\hat{h}_f(P) - h(P)| \leq C$ for all $P \in \P^1(k)$. 
\end{enumerate}
In our geometric setting of $k = \C(X)$, from equation (\ref{degree}) the canonical height of the point $P$ is nothing more than the degree growth of the iterates,
\begin{equation} \label{degree growth}
	\hat{h}_f(P) = \lim_{n\to\infty} \frac{\deg(f^n(P))}{d^n}
\end{equation}
where each iterate $f^n(P)$ is viewed as a map $X \to \P^1$, defined by $t\mapsto f_t^n(P(t))$.   

Observe from the definition of $\hat{h}_f$ that a preperiodic point $P$ will always satisfy $\hat{h}_f(P) = 0$; this is the implication (2) $\implies$ (3) in Theorem \ref{mine}.  As a consequence, property (H2) immediately implies the Proposition \ref{finite degree} stated in Lecture 1.  

The converse statement, that height 0 implies the point is preperiodic is the content of (3) $\implies$ (2) in Theorem \ref{mine}; it is false without the hypothesis of non-isotriviality (though it is always true when working over number fields \cite[Corollary 1.1.1]{Call:Silverman}).  

\subsection{Idea of the proof of Theorem \ref{mine}}
As we have just stated, the implication (2) $\implies$ (3) follows from the definition of the canonical height.  
The implication (3) $\implies$ (2) and a proof of the finiteness of all rational preperiodic points can be obtained with the same argument used in Lecture 1 \S\ref{finiteness}.  In rereading that proof, we see that the structure of the elliptic curve is not used, and all of the statements hold for an arbitrary non-isotrivial rational function $f \in k(z)$.  To obtain the result under the weaker hypothesis of non-isotriviality over $k$, one needs an extra (easy) argument about automorphisms; see \cite[Lemma 2.1]{D:stableheight}.  

Condition (2) implies condition (1), immediately from the definitions.  So it remains to show that condition (1) implies the others.  This is the piece of the theorem generalizing Theorems \ref{McMullen} and \ref{DF critical}.  Because of the arguments in \S\ref{finiteness}, it suffices to show that stability of $(f,P)$ implies that the degrees of the iterates $f^n(P)$ are bounded.  We proceed by contradiction:  assume there is a sequence of iterates $f^{n_k}(P)$, $k\geq 1$, with unbounded degrees as maps $X \to \P^1$.  The stability on the punctured Riemann surface $V\subset X$ guarantees the existence of a further subsequence, which we will also denote $f^{n_k}(P)$, that converges locally uniformly on $V$ to a holomorphic map 
	$$\phi: V \to \P^1.$$
The graph of the function $P: X\to \P^1$ determines an algebraic curve $\Gamma_P$ inside the complex surface $X\times \P^1$.  Let $\Gamma_k$ be the graph of $f^{n_k}(P)$.  Then $\Gamma_k$ will converge (over compact subsets of $V$) to the graph $\Gamma_\phi$ of $\phi$.  If $\deg f^{n_k}(P) \to \infty$ as $k\to\infty$, then these graphs $\Gamma_k$ will be fluctuating wildly near the punctures of $V$ for all $k>>0$.  

To complete the proof that (1) implies (2), we now exploit the algebraicity of $f$:  note that $f$ induces a rational map
	$$F:  X\times \P^1 \dashrightarrow X\times \P^1$$
given by $(t, z) \mapsto (t, f_t(z))$, with indeterminacy points over the punctures of $V$.  The ``depth" of the indeterminacy of $F$ is controlled by the order of vanishing of the resultant of $f_t$ at each of these punctures; a precise statement appears as \cite[Lemma 3.2]{D:stableheight}.  In particular, the iterates of $F$ cannot be too wild near these punctures; a precise estimate on this control appears as \cite[Proposition 3.1]{D:stableheight}, formulated in terms of the homogeneous escape-rate function shown in equation (\ref{escape rate}).  It follows that the graphs $\Gamma_k$ cannot be fluctuating near the punctures of $V$, contradicting the assumption of unbounded degree.

\bigskip
\section{Dynamical moduli spaces and unlikely intersections}

In this final lecture, we discuss some open questions about the dynamical systems $f: \P^1\to\P^1$ inspired by numerous results and questions in arithmetic and algebraic geometry.  The general question in its vague, heuristic form is:  
\begin{quote}  Suppose $M$ is a complex algebraic variety, perhaps a moduli space of objects.  Suppose $Z \subset M$ is a subset of ``special" points within $M$.  If a complex algebraic curve $C\subset M$ passes through infinitely many points of $Z$, then what is special about $C$?  
\end{quote}
One can, and does, ask the same question for algebraic subvarieties $V\subset M$ of arbitrary dimension, where ``infinitely many" is replaced by ``a Zariski dense set of" -- but for this lecture, we will concentrate on the 1-dimensional subvarieties. 

A classical example is the Manin-Mumford Conjecture, which is a theorem of Raynaud \cite{Raynaud:1, Raynaud:2}.  Let $M$ be an abelian variety, a projective complex torus, and let $Z$ be the set of torsion points in $M$.  Then $C\cap Z$ is infinite if and only if $C$ is either an abelian subvariety or a torsion-translate of an abelian subvariety.  In this example, the set $Z$ of torsion points is dense in $M$.  (The original formulation of the question in this setting was about curves $C$ of genus $>1$ embedded in their Jacobians.)

A modular example was proved by Andr\'e \cite{Andre:finitude}.  Suppose $M = \C^2$ parameterizes pairs of elliptic curves $(E_1, E_2)$ by their $j$-invariants.  Recall that the $j$-invariant of an elliptic curve uniquely determines the curve up to isomorphism (over $\C$).  Let $Z$ be the subset of pairs $(E_1, E_2)$ where both $E_i$ have complex multiplication (CM).  This means that $E_i$ has endomorphisms in addition to the usual $P \mapsto n\cdot P$ with $n\in \Z$.  For example, the square torus $\C/(\Z\oplus i\Z)$ has a complex automorphism given by $z\mapsto i z$, and so it is CM.  The rectangular torus $\C/(\Z\oplus i\sqrt{2} \Z)$ has an endomorphism of degree 2 given by $z \mapsto i\sqrt{2} z$, so it is also CM.  For this example, again $Z$ is a dense subset of $M = \C^2$, though the elements of $Z$ are arithmetically ``special"; see \cite{Cox:primes} for an accessible introduction to the rich theory of CM elliptic curves.  

Andr\'e's theorem states that a complex algebraic curve $C\subset M = \C^2$ has infinite intersection with $Z$, the set of CM pairs in $\C^2$, if and only if $C$ is either 
\begin{enumerate}
\item a vertical line $\{E_0\}\times \C$ where $E_0$ has complex multiplication, 
\item a horizontal line $\C\times \{E_0\}$ where $E_0$ has complex multiplication, or 
\item the modular curve $Y_0(N)$, consisting of pairs $(E_1, E_2)$ for which there exists an isogeny (a covering map, which is also a group homomorphism) of degree $N$, for any positive integer $N$.
\end{enumerate}

Generalizations of Andr\'e's theorem and related results led to the development of the Andr\'e-Oort Conjecture, characterizing the ``special subvarieties" of Shimura varieties; see results of Pila and Tsimerman \cite{Pila:AO, Tsimerman:AO}.  Further conjectures were put forward by Pink and Zilber; see also \cite{Zannier:book}.  In the study of dynamical systems $f: \P^1\to \P^1$, there are certain ``special" maps that play the role of elliptic curves with complex multiplication in their moduli space, the postcritically finite maps described below.

\subsection{A dynamical Andr\'e-Oort question}
Let $\M_d$ denote the moduli space of rational maps 
	$$f: \P^1\to\P^1$$
of degree $d\geq 2$; by definition, it is the set of all conformal conjugacy classes of rational functions.  It is the quotient of $\Rat_d$ (defined in \S\ref{definitions}) by the group $\PSL_2\C$ of M\"obius transformations acting by conjugation.  The space $\M_d$ is a complex affine algebraic variety of dimension $2d-2$; see  \cite{Silverman:moduli} for more information about its structure.  In the case of $d=2$, it is isomorphic to $\C^2$ \cite{Milnor:quad}.  

A map $f: \P^1 \to \P^1$ of degree $d$ is said to be {\em postcritically finite} if each of its $2d-2$ critical points has a finite forward orbit.  The set of postcritically finite maps within $\M_d$ is not dense (in the usual topology), but it is Zariski dense, meaning that it does not lie in a complex algebraic hypersurface in $\M_d$.  See Theorem A of the Appendix for a proof.  

\begin{question}  \label{dynamical AO}
Suppose $\M_d$ is the moduli space of complex rational maps $f: \P^1\to \P^1$ of degree $d\geq 2$.  Let $\mathrm{PCF}\subset\M_d$ be the set of postcritically-finite maps of degree $d$.  Which algebraic curves $C\subset \M_d$ pass through infinitely many elements of $\mathrm{PCF}$?
\end{question}

\noindent 
I first heard this question asked by Bjorn Poonen, during the series of lectures by Joe Silverman in May 2010 that led to the writing of \cite{Silverman:moduli}.  At the time, I put forward a roughly-stated conjectural answer:  any such $C$ must be defined by ``critical orbit relations."  This question and the conjecture are discussed briefly in \cite{Silverman:moduli} and at greater length in \cite{BD:polyPCF} and \cite{D:stableheight}.

To formulate a precise conjectural answer to Question \ref{dynamical AO}, it helps to work in the branched cover $\M_d^{cm} \to \M_d$, consisting of conformal equivalence classes of critically-marked rational maps $(f, c_1, \ldots, c_{2d-2})$, where each $c_i$ is a critical point of $f$.  An example of a critical orbit relation is an equation of the form 
\begin{equation} \label{orbit relation}
	f^n(c_i) = f^m(c_j).
\end{equation} 
for a pair of non-negative integers $(n,m)$ and any $i, j \in \{1, \ldots, 2d-2\}$.  To answer Question \ref{dynamical AO}, it turns out we will need to consider a larger class of critical orbit relations, allowing for symmetries of the map $f$.  See \cite[Theorem 1.2]{BD:polyPCF} and the examples presented there; see also the discussion in \cite[\S6.2]{D:stableheight} and part (3) of Conjecture \ref{UI} below.  

\begin{remark}
There is an interesting connection between critical orbit relations of the form \eqref{orbit relation} and stability in holomorphic families $f: X\times\P^1\to \P^1$ \cite[Theorem 2.7]{McS:QCIII}.  McMullen and Sullivan proved that the family is structurally stable -- a stronger notion of stability than that introduced in \S\ref{stable}, requiring topological conjugacy on all of $\P^1$ -- if and only if any critical orbit relation that holds at a parameter $t_0 \in X$ persists throughout $X$.  
\end{remark}

The curves $C$ in $\M_d$ containing infinitely many elements of $\mathrm{PCF}$ will be called {\em special}.  Two explicit examples of special curves are provided by the examples from Lecture 1.  First, consider the subspace of polynomials within the space of quadratic maps $\M_2$.  There are infinitely many postcritically finite polynomials in the family $f_t(z) = z^2+t$, consisting of all solutions $t$ to equations of the form 
	$$f_t^n(0) = f_t^m(0)$$
with $n> m \geq 0$.  The family of quadratic polynomials is itself defined by a critical orbit condition in the moduli space $\M_2$:  that there exists a critical point $c$ for which $f(c)=c$.  (This is the critical point that lies at infinity for the polynomial.)  A second example of a special curve is given by the family of Latt\`es maps (\ref{Lattes family}) in degree 4 (or in any square degree), since every Latt\`es map is postcritically finite.  Their critical orbit relations, where under two iterates all critical points have landed on a fixed point, characterize the Latt\`es maps in $\M_4$ \cite{Milnor:Lattes}, so this family is also defined by critical orbit relations.  

To further illustrate the conjectural answer to Question \ref{dynamical AO}, again consider the moduli space of quadratic rational maps $\M_2 \iso \C^2$.  In \cite{Milnor:quad}, Milnor introduced the family of curves $\Per_1(\lambda) \subset \M_2$, $\lambda\in\C$, consisting of all maps $f$ with a fixed point of multiplier $\lambda$; that is, maps for which there exists a point $p$ with $f(p) = p$ and $f'(p)=\lambda$.  Note that $\Per_1(0)$ is defined by a critical orbit relation; it is precisely the family of quadratic polynomials.  These curves sweep out all of $\M_2$ as $\lambda$ varies.  Since $\mathrm{PCF}\subset \M_2$ is Zariski dense, we know that infinitely many of these lines $\Per_1(\lambda)$ must contain points of $Z$.  Nevertheless:

\begin{theorem} \label{QPer1} \cite{DWY:QPer1}
The curve $\Per_1(\lambda)\subset \M_2$ contains infinitely many postcritically-finite maps if and only if $\lambda=0$.
\end{theorem}

A second illustrative result is closely related to Question \ref{dynamical AO}, though not exactly a particular case, paralleling the result of Andr\'e (stated above) for pairs of elliptic curves with complex multiplication.

\begin{theorem} \label{Holly} \cite{GKNY}
Let $M = \C^2$ parameterize pairs of quadratic polynomials $(z^2+a, z^2+b)$ by $(a,b) \in \C^2$.  Let $Z$ be the set of all postcritically-finite pairs.  Then an algebraic curve $C\subset M$ has infinite intersection with $Z$ if and only if 
\begin{enumerate}
\item $C$ is a horizontal line $\C\times\{c_0\}$ where $z^2+c_0$ is postcritically finite; 
\item $C$ is a vertical line $\{c_0\}\times\C$ where $z^2+c_0$ is postcritically finite; or
\item $C$ is the diagonal $\{a=b\}$.
\end{enumerate}
\end{theorem}

\noindent
Note that in the setting of Theorem \ref{Holly}, there are no non-trivial symmetries among the quadratic polynomials, so there is no infinite family of ``modular curves" as in the case Andr\'e studied.

\begin{remark}   \label{progress}
The first result answering a special case of Question \ref{dynamical AO} appears in \cite{Ghioca:Hsia:Tucker}, treating certain families of polynomials.  More families of polynomials were examined in \cite{BD:polyPCF} where a general conjectural answer was proposed, giving a full treatment of symmetries and critical orbit relations.  Most recently, about a year after these lectures were given, a complete answer to Question \ref{dynamical AO} was obtained in the case of curves $C$ in $\MPoly_3$, the moduli space of cubic polynomials \cite{Favre:Gauthier:cubics, Ghioca:Ye:cubics}.  There are extra tools available to study the dynamics of polynomials, while the case of curves $C$ in $\M_d$ is mostly wide open.

The proofs of Theorems \ref{QPer1} and \ref{Holly}, and the results just mentioned, all follow the same general strategy, though the technical details and difficulties differ.  An outline of these proofs is presented below in \S\ref{outline}.
\end{remark}

\subsection{Unlikely intersections in dynamical moduli spaces}  \label{unlikely}
Question \ref{dynamical AO} is focused on the orbits of critical points, because of the dynamical importance of critical orbit behavior, but we could also study orbits of arbitrary points.  In the context of elliptic curves, such a study has been carried out; the following result of Masser and Zannier addressed a special case of the Pink and Zilber conjectures.

\begin{theorem}  \cite{Masser:Zannier}  \label{MZ}
Let $E_t = \{y^2 = x(x-1)(x-t)\}$ be the Legendre family of elliptic curves, $t\in \C\setminus\{0,1\}$, and let $P$, $Q$ be points on $E$ with $x$-coordinates equal to $a,b \in \C(t)$, respectively, with $a,b \not= 0,1,\infty,t$.  If the set $\{t:  P_t, Q_t \mbox{ are both torsion on }E_t\}$ is infinite, then there exist nonzero integers $n,m$ so that $nP+mQ = 0$ on $E$.  
\end{theorem}

\noindent
The articles \cite{Masser:Zannier:CR, Masser:Zannier, Masser:Zannier:2} discuss the importance of Theorem \ref{MZ} and several generalizations and implications.  Zannier's book  \cite{Zannier:book} provides an overview.  

From Lecture 1, we will recall that torsion points on elliptic curves can be studied dynamically.  The article \cite{DWY:Lattes} contains a dynamical proof of Theorem \ref{MZ}, analyzing features of the Latt\`es family (\ref{Lattes family}), building upon the computations of \S\ref{Lattes computation}.  

More generally, suppose $k= \C(X)$ for a compact Riemann surface $X$ and fix a non-isotrivial $f\in k(z)$.  Fix points $P, Q\in \P^1(k)$ and define sets
	$$S_{f,P} = \{t\in X:  P_t  \mbox{ has finite forward orbit under } f_t\}$$
and 
	$$S_{f,Q} = \{t\in X:  Q_t  \mbox{ has finite forward orbit under } f_t\}.$$
Theorem \ref{mine}, combined with a standard complex-dynamic argument using Montel's Theorem on normal families, shows that the sets $S_{f,P}$ and $S_{f,Q}$ are infinite in $X$.  Indeed, if the point $P$ is persistently preperiodic for $f$ then it has finite orbit for all $t$.  Otherwise, by Theorem \ref{mine}, the bifurcation locus $\mathcal{B}(f,P)$ must be nonempty in $X$.  Then, by Montel's theory of normal families, in every open set intersecting $\mathcal{B}(f,P)$, there must be parameters $t$ for which an iterate $f_t^n(P(t))$ lands on a preimage of $P(t)$, making the point $P(t)$ periodic.  Similarly for $Q$.

If the points $P$ and $Q$ are ``independent" in some dynamical sense, then we should expect the intersection of $S_{f,P}$ and $S_{f,Q}$ to be small, or perhaps empty.  Inspired by Theorem \ref{MZ}, we ask:

\begin{question} \label{dynamical UI}
Let $k = \C(X)$, and fix a non-isotrivial $f\in k(z)$ and $P, Q\in \P^1(k)$.  If the set
	$$S_{f,P}\cap S_{f,Q} = \{t\in X:  P_t \mbox{ and } Q_t \mbox{ both have finite orbit under } f_t\}$$
is infinite, then what can we say about the triple $(f, P,Q)$?
\end{question}

I suspect the following might be true:

\begin{conjecture} \label{UI}
Let $k = \C(X)$, and fix non-isotrivial $f\in k(z)$ and $P, Q\in \P^1(k)$.  Suppose that neither $P$ nor $Q$ is persistently preperiodic (having finite orbit for every $t \in X$).  Then the following are equivalent:
\begin{enumerate}
\item	$|S_{f,P} \cap S_{f,Q}| = \infty$, 
\item	$S_{f,P} = S_{f,Q}$
\item there exist $A, B\in \kbar(z)$ and integer $\ell \geq 1$ so that 
	$$f^\ell\circ A = A \circ f^\ell, \; \; f^\ell\circ B = B \circ f^\ell, \; \; \mbox{ and } \; A(P) = B(Q)$$ 
\end{enumerate}
\end{conjecture}

\noindent
The conclusion of Conjecture \ref{UI} (3) is that {\em $P$ and $Q$ have the same grand orbit under $f$, up to symmetries}.  For example, $A$ and $B$ might be equal to iterates of $f$.  Compare \cite[Theorem 1.2]{BD:polyPCF}.  It is known that (3) implies (2) implies (1).  Details appear in \cite{D:stableheight}.

Note that Theorem \ref{MZ} falls under the umbrella of Conjecture \ref{UI}:  we fix $f$ to be our flexible Latt\`es family (\ref{Lattes family}) of degree 4.  The hypothesis that the $x$-coordinates of $P$ and $Q$ are not $0,1,t,\infty$ guarantees that $P$ and $Q$ are not torsion on $E$ (for all $t$), as we can see from the discussion in \S\ref{Lattes computation}.  The relation $nP + mQ = 0$ on $E$ corresponds to an orbit relation on the $x$-coordinates of $P$ and $Q$, under the symmetries $A = f_{[n]}$ and $B = f_{[-m]}$ (induced from multiplication by $n$ and by $-m$ on $E$) with $\ell =1$.

\subsection{Proof strategy}  \label{outline}
The proofs of Theorems \ref{QPer1} and \ref{Holly}, and the dynamical proof of Theorem \ref{MZ} follow the same general strategy, built on the mechanism developed in \cite{BD:preperiodic, Ghioca:Hsia:Tucker, GHT:rational, BD:polyPCF}.  A possible strategy for proving Conjecture \ref{UI} might involve the following three steps:
\begin{enumerate}
\item[Step 1.]  The application of an arithmetic dynamical equidistribution theorem on the Berkovich analytification of $X$, as in \cite{Baker:Rumely:equidistribution, FRL:equidistribution, ChambertLoir, Yuan:equidistribution, Thuillier:these}.  The elements of $S_{f,P}$ are expected to be uniformly distributed in $X$ with respect to the bifurcation measure $T_P$ (defined in \S\ref{point bifurcations}).  Moreover, as $f$ and $P$ are defined over $K(X)$ for a finitely-generated field over $\Q$, any $\Gal(\Kbar/K)$-invariant infinite subset of $S_{f,P}$ is also expected to be equidistributed with respect to $T_P$;  in particular, the hypothesis that $|S_{f,P}\cap S_{f,Q}| = \infty$ should guarantee that $T_P = T_Q$ (and in turn, that $S_{f,P} = S_{f,Q}$ or that their symmetric difference is finite).  
\item[Step 2.]  Deduce from $S_{f,P} = S_{f,Q}$ or from $T_P =T_Q$ that there is a dynamical relation between $P$ and $Q$.  A general problem in the analytic setting is posed as Question \ref{critical question}.  For example, in the case of polynomial families, the existence of the B\"ottcher coordinate near $\infty$ was useful in producing the desired relations in \cite{BD:preperiodic, BD:polyPCF}.
\item[Step 3.]  A characterization of symmetries of $f$ to obtain the exact form of the relation obtained in Step 2.  For example, in \cite{GKNY}, the authors proved a rigidity result about the Mandelbrot set to obtain the complete list of special curves for Theorem \ref{Holly}.  In \cite{BD:polyPCF}, we appealed to the results of \cite{Medvedev:Scanlon} which built upon the work of Ritt in \cite{Ritt:decompositions}, to obtain a simple form for the critical orbit relations of polynomials.  
\end{enumerate}

\bigskip
\section*{Appendix.  Zariski density of hyperbolic postcritically finite maps} \label{density}

Let $\M_d$ denote the moduli space of complex rational maps 
	$$f: \P^1\to \P^1$$
of degree $d\geq 2$.  In this Appendix, we provide a proof that the set PCF of postcritically finite maps in the moduli space $\M_d$ is Zariski dense.  Recall that $f$ is postcritically finite if each of its critical points has a finite forward orbit.  In fact, we prove a stronger statement, showing that the following subset of the postcritically finite maps is also Zariski dense.  Let 
	$$\mathrm{HPCF} = \{f \in \M_d: \mbox{ every critical point of } f \mbox{ is periodic}\}.$$
The H in HPCF stands for hyperbolic, since all such maps are expanding on their Julia sets; see, e.g., \cite[\S3.4]{McMullen:CDR} for definitions and the characterizations of hyperbolicity for rational maps $f: \P^1\to \P^1$.

A similar proof to the one given below appears in \cite[Proposition 2.6]{BD:polyPCF}, proving that the set PCF is Zariski dense in the moduli space of polynomials; the Zariski-density of PCF is also a consequence of \cite[Theorem 1.6]{D:stableheight}.  In fact, the same key ingredients appear also in the proofs of \cite[Corollary 5.3]{Dujardin:supports}, \cite[Proposition 3.7]{Gauthier:strong}, and \cite[Lemma 4.2]{Gauthier:higher}.  The proof given here is not the only approach to showing Zariski density.  The sketch of an alternative argument appeared in \cite[Proposition 6.18]{Silverman:moduli}, building on Epstein's transversality theory.  Or, one may appeal to the fact that the bifurcation current has continuous potentials (Theorem \ref{current theorem} of Lecture 2), so the bifurcation measure (the top wedge power of the bifurcation current) does not assign positive mass to algebraic subvarities \cite[Proposition 117]{Berteloot:lectures}.  Combined with the fact that the (hyperbolic) postcritically finite maps accumulate everywhere in the support of this measure, we obtain Zariski density; see, e.g., \cite[Theorem 134]{Berteloot:lectures} or \cite[Theorem 3.2]{Dujardin:bifurcation}.

I believe the strategy here is the most direct, relying on the least amount of additional machinery, at least to prove the weaker statement that PCF is dense in $\M_d$.  The ideas are in line with the concepts introduced in these lectures, relying only on \cite[Theorem 2.5]{Dujardin:Favre:critical} (or its generalized form, Theorem \ref{mine} of Lecture 3).  To obtain the density of HPCF, I use one additional nontrivial ingredient in the proof below, namely Thurston Rigidity \cite{Douady:Hubbard:Thurston}. 

\begin{theoremA}  \label{Zariski}
For every $d\geq 2$, the set $\mathrm{HPCF}$ is Zariski dense in $\M_d$.  
\end{theoremA}

\proof
Let $S$ be any proper algebraic subvariety of $\M_d$, and let $\Lambda$ be its complement.  It suffices to show that there exists an element of $\mathrm{HPCF}$ in $\Lambda$.  

Since $\M_d$ is an irreducible complex affine variety, we see that $\Lambda$ is an irreducible quasi-projective algebraic variety.  Choose an algebraic family of rational maps $f: V\times\P^1\to \P^1$ that projects to $\Lambda$ in $\M_d$, where $V$ is also an irreducible quasi-projective complex algebraic variety; for example, $V$ might be the preimage of $\Lambda$ in the space $\Rat_d$ of all rational functions of degree $d$.  Passing to a branched cover of $V$ if necessary, we may assume that all $2d-2$ critical points of $f$ can be holomorphically parameterized, by $c_1, \ldots, c_{2d-2}:  V\to \P^1$.  Note that $\dim_\C V \geq 2d-2 = \dim_\C \M_d$.  

Consider the critical point $c_1$.  Observe that $c_1$ cannot be persistently preperiodic on all of $V$, because it would then be persistently preperiodic for every $f \in \M_d$; this is absurd, as there exist maps in every degree for which all critical points have infinite orbit.  Thus, by \cite[Theorem 2.5]{Dujardin:Favre:critical}, the pair $(f, c_1)$ must be undergoing bifurcations.  Consequently, there will be a parameter $v_1\in V$ at which $c_1(v_1)$ is periodic for $f_{v_1}$; this conclusion requires an application of Montel's Theorem, as in the proof of \cite[Proposition 2.4]{Dujardin:Favre:critical} (or the discussion above, just before Question \ref{dynamical UI}).  

Suppose $v_1$ satisfies the equation $f^{n_1}_{v_1}(c_1(v_1)) = c_1(v_1)$.  Let $V_1 \subset V$ be an irreducible component of the subvariety defined by this equation, containing $v_1$.  Then $V_1$ is a quasiprojective variety of codimension $1$ in $V$, and the restricted family $f: V_1\times\P^1\to \P^1$ is an algebraic family of rational maps of degree $d$, with marked critical points, for which $c_1$ is persistently periodic.  Moreover, if we denote by $\Lambda_1$ the projection of $V_1$ in $\M_d$, then $\Lambda_1$ will have codimension 1 in $\Lambda$ (because not all $f$ in $\M_d$ have a periodic critical point, so the equation $f^{n_1}(c_1) = c_1$ defines a hypersurface).

We continue inductively.  Suppose $V_k$ is a quasiprojective subvariety of codimension $\leq k$ in $V$, and $f: V_k\times\P^1\to \P^1$ is an algebraic family of rational maps (with marked critical points) for which the critical points $c_1, \ldots, c_k$ are persistently preperiodic.  We also assume that the projection $\Lambda_k$ of $V_k$ in the moduli space $\M_d$ has the same codimension in $\Lambda$ as that of $V_k$ in $V$.   If $c_{k+1}$ is persistently preperiodic along $V_k$, set $V_{k+1} = V_k$.  If not, we know from \cite[Theorem 2.5]{Dujardin:Favre:critical} that $(f,c_{k+1})$ must be bifurcating along $V_k$, and we can guarantee the existence of a parameter $v_{k+1} \in V_k$ at which $c_{k+1}$ is periodic.  We let $V_{k+1} \subset V_k$ be an irreducible component of the subvariety defined by the critical orbit relation on $c_{k+1}$, containing the parameter $v_{k+1}$.   Then $V_{k+1}$ has codimension 1 in $V_k$.  It follows that $\Lambda_{k+1}$ has codimension 1 in $\Lambda_k$, because critical orbit relations are constant along the fibers of the projection $V\to \M_d$.  

Now, the induction argument has produced for us an algebraic family 
	$$f:  V_{2d-2} \times\P^1\to \P^1$$
for which all critical points are persistently preperiodic and $c_1$ is persistently periodic.  The variety $V_{2d-2}$ projects to a quasiprojective variety $\Lambda_{2d-2} \subset \mathrm{PCF} \subset \M_d$.   From the construction, the codimension of $\Lambda_{2d-2}$ in $\M_d$ is no larger than the number of periodic critical points of any $f\in \Lambda_{2d-2}$.  In particular, the codimension of $\Lambda_{2d-2}$ is $\leq 2d-2$, so the variety $\Lambda_{2d-2}$ is nonempty.  This shows that there must be a postcritically finite map in $\Lambda$.  

This completes the proof that PCF is Zariski dense in $\M_d$.

If all critical points for $f\in \Lambda_{2d-2}$ are periodic, we have shown that HPCF is also Zariski dense in $\M_d$, since then $\Lambda_{2d-2}$ is a nonempty subset of HPCF in $\Lambda$.  

Suppose that not all critical points are periodic for the maps $f\in \Lambda_{2d-2}$.  Then, by the construction, the codimension of $\Lambda_{2d-2}$ in $\M_d$ must be $< 2d-2$, so $\Lambda_{2d-2}$ defines a positive-dimensional family of postcritically finite maps in $\M_d$ for which at least one critical point (namely, $c_1$) is periodic.  But this is impossible, due to Thurston's Rigidity Theorem, which states that the flexible Latt\`es maps are the only positive-dimensional family of postcritically finite maps in $\M_d$, and all critical points for any Latt\`es map are strictly preperiodic.  We conclude, therefore, that $\Lambda_{2d-2}$ has codimension $2d-2$ and that $\Lambda_{2d-2} \subset \mathrm{HPCF}$.  This concludes our proof.  
\qed

\bigskip \bigskip
\def\cprime{$'$}

\bigskip\bigskip

\end{document}